\documentclass[letterpaper, 10 pt, conference]{ieeeconf}  % Comment this line out if you need a4paper

\IEEEoverridecommandlockouts                              % This command is only needed if 
\usepackage{graphics} % for pdf, bitmapped graphics files
\usepackage{epsfig} % for postscript graphics files
\usepackage{mathptmx} % assumes new font selection scheme installed
\usepackage{amsmath} % assumes amsmath package installed
\usepackage{amssymb}  % assumes amsmath package installed

\title{\LARGE \bf
A Unified Approach to Configuration-based Dynamic Analysis of Quadcopters for Optimal Stability
}

\author{Mojtaba Hedayatpour$^{1}$, Mehran Mehrandezh$^{1}$ and Farrokh Janabi-Sharifi$^{2}$% <-this % stops a space
\thanks{$^{1}$M. Hedayatpour and M. Mehrandezh are with University of Regina, Regina, SK S4S 0A2 Canada. 
        {\tt\small \{{hedayatm,mehran.mehrandezh\}}@uregina.ca}}%
\thanks{$^{2}$F. Janabi-Sharifi is with Ryerson University, Toronto, ON M5B 2K3 Canada. 
        {\tt\small fsharifi@ryerson.ca}}%
}

\begin{document}

\maketitle
\thispagestyle{empty}
\pagestyle{empty}

%%%%%%%%%%%%%%%%%%%%%%%%%%%%%%%%%%%%%%%%%%%%%%%%%%%%%%%%%%%%%%%%%%%%%%%%%%%%%%%%
\begin{abstract}

A special type of rotary-wing Unmanned Aerial Vehicles (UAV), called Quadcopter have prevailed to the civilian use for the past decade. They have gained significant amount of attention within the UAV community for their redundancy and ease of control, despite the fact that they fall under an under-actuated system category. They come in a variety of configurations. The ``+" and ``x" configurations were introduced first. Literature pertinent to these two configurations is vast. However, in this paper, we define 6 additional possible configurations for a Quadcopter that can be built under either ``+" or ``x" setup. These configurations can be achieved by changing the angle that the axis of rotation for rotors make with the main body, i.e., fuselage. This would also change the location of the COM with respect to the propellers which can add to the overall stability. A comprehensive dynamic model for all these configurations is developed for the first time. The overall  stability for these configurations are addressed. In particular, it is shown that one configuration can lead to the most statically-stable platform by adopting damping motion in Roll/Pitch/Yaw, which is described for the first time to the best of our knowledge.

\end{abstract}

%%%%%%%%%%%%%%%%%%%%%%%%%%%%%%%%%%%%%%%%%%%%%%%%%%%%%%%%%%%%%%%%%%%%%%%%%%%%%%%%
\section{INTRODUCTION}

Multi-copter unmanned aerial vehicles (UAVs) with vertical take-off and landing (VTOL) capability  are becoming more popular due to the ease of their  operation. They come in a variety of shapes and configurations. Among them, quadcopters are gaining a lot of attention due to their simple structure and ease of control [1], quadcopters are used in application domains such as: aerial photogrammetry, aerial inspection of infrastructure, precision agriculture, immersive televising of sports events, and object delivery [1]-[3].

They usually come in two configurations, namely ``+" and ``x" configurations [4]. The main advantage of the ``x" configuration is mainly due to its open frontal area that facilitates for employment of occlusion-free forward-looking imaging sensors. Although, their dynamics would be different, not much attention has been given to their subtle differences within the research community.

Quadcopters with fixed rotors fall under the under-actuated and non-holonomic flying machine categories. Adoption of a larger number of rotors and/or adding the tilting effect on them for on-the-fly thrust vectoring can lead to fully-actuated holonomic machines at the cost of making them mechanically more complicated and less power efficient. 

There have been some studies on: (i) building UAVs using variable-pitch blades [5]; (ii) configuring rotors to yield non-parallel thrust vectors [6]-[8] and (iii) designing multi-copter UAVs with rotors that can tilt on the fly [9] (iv) building multi-copter UAVs with rotors fixedly mounted with an angle with respect to the fuselage [10]. However, very little attention has been given to calculating the optimal configuration in quadcopters with fixed rotors for highest static and dynamic stability. In this paper, we attempt to look at all possible controllable configurations for a quadcopter with fixed rotors and analyze their stability attributes in a quantitative fashion for the first time. We also provide a unified dynamic model for all the possible configurations from which special cases can be deducted.

Literature pertinent to the mathematical modeling of quadcopters and their flight control is vast, [2]-[8].  In our derivation, we assume a full model of the gyroscopic moments for the first time. More specifically, we derive the dynamic model of quadcopters assuming that: (A1) the thrust vector for each rotor would make a non-zero angle with the vertical axis (i.e., the sagittal suture) of the quadcopter; and (A2) the center-of-mass (COM) of the quadcopter does not lie on the same plane where the center-of-mass of all motors lie on (blue plane shown in Fig. 1). However, we still assume that the quadcopter under study has two axes of congruency (see Fig. 1).

%%thpb
   \begin{figure}[b]
      \centering
      \includegraphics[scale=0.8]{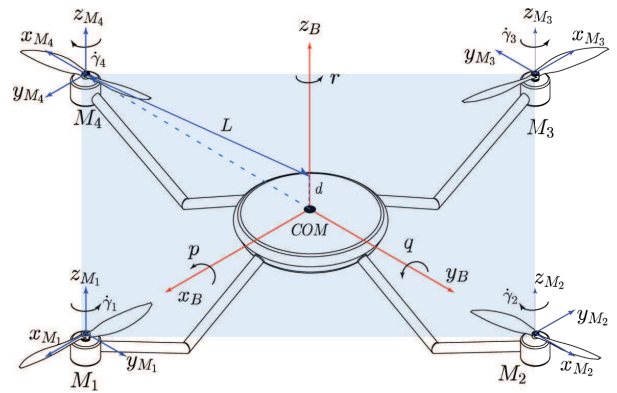}
      \caption{Quadcopter in ``+" configuration. Body frame is shown in blue and is attached to the center-of-mass of the quadcopter. A frame, shown in blue, is attached to each motor in order to determine orientation of the motors with respect to body frame. Motors are located at distance L and d from z-axis and x-y plane of the body frame respectively.}
      \label{quad}
   \end{figure}

The angle between the thrust vector of each rotor and the vertical axis of the fuselage is further divided into: (1) the dihedral angle, and (2) twist (i.e., lateral tilting) angle ( Figs. 2 and 3). We assume that the central hub of all four blades lie on a flat horizontal plane (blue plane in Fig. 1), called flat plane from this point on, from which the location of the COM is referenced (i.e., the COM can be either above, below, or right on this plane).

The dynamic model developed in this paper will, therefore, have three additional terms in comparison to that in the flat quadcopters (this is the term used for the original quadcopters, where the COM and the rotor hubs are all on the same plane), as: dihedral angle $\beta_i$, twist angle $\alpha_i$, and the distance between the COM and the flat plane $d$ (please note that $d$ could take positive and negative values, measured in z-direction of the body frame). In existing flat-model of quadcopters one has: $\beta_i=\alpha_i=d=0$. 

We will show that the flat model of quadcopters does not render itself as the most statically and dynamically-stable configuration. For instance, by adding a dihedral angle to the blades' thrust vectors, one can achieve better rolling stability in forward flight. Also, the twist angles in the blades would yield faster yaw dynamics without compromising the overall stability of the system. Furthermore, a positive value of $d$ (i.e., positioning the COM of the quadcopter below the flat plane), one can achieve an open-loop roll/pitch stable configuration. 

We use Newton's method for driving the dynamic model of the quadcopter. Also, without the loss of generality, we assume a ``+" configuration. The rest of the paper is organized as follows: In section II, the derivation of equations of motion is presented. In section III, the effects of having dihedral and twist angles are given. In section IV, stability analysis for six different configurations is provided and compared with that in the flat-model quadcopters. Conclusions and future work are presented in Section V.

   \begin{figure}[b]
      \centering
      \includegraphics[scale=0.9]{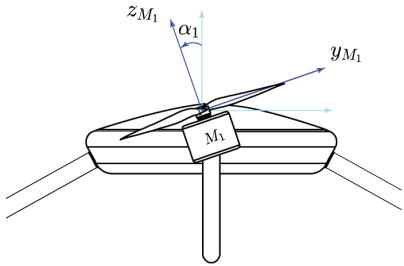}
      \caption{Twist angle $\alpha_1$ about the x-axis of the motor frame $M_1$.}
      \label{alphai}
   \end{figure}
   
   \begin{figure}[b]
      \centering
      \includegraphics[scale=0.6]{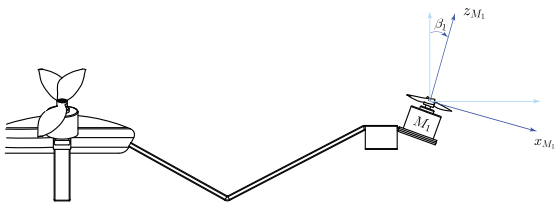}
      \caption{Dihedral angle $\beta_1$ about the y-axis of the motor frame $M_1$.}
      \label{betai1}
   \end{figure}

\section{EQUATIONS OF MOTION}

\subsection{Notation \& Parameters}

Since there are many rotation matrices involved in this modeling, straight boldface letter \textbf{R} is only reserved for rotation matrices. The rotation from frame A to frame B is expressed as $^B$\textbf{R}$_A$. Also $^B{\pmb{\omega}}_{P_i,I}$ indicates that $\pmb{\omega}$ belongs to the $i^{th}$ propeller with respect to an inertial frame $I$ and is expressed in the body frame $B$. \textbf{R}$_A(\theta)$ , represents a rotation matrix about axis $A$ by angle $\theta$.

\subsection{Frames \& Transformations}

The body frame $^BO- {^Bx ^By ^Bz}$ (red color in Fig. 1) is attached to the center of mass of the vehicle. Four frames named $^{M_i}O- {^{M_i}x ^{M_i}y ^{M_i}z}$ (blue color in Fig. 1) are attached to motors. Motors are turning with angular velocities $\dot{\gamma_i}$ $(i=1,2,...,4)$ about axis $z_{M_i}$. Position of the vehicle is expressed in the inertial frame $I$.

Orientation of the body frame with respect to the inertial frame can be captured by the rotation matrix from body frame to inertial frame $^I$\textbf{R}$_B$ . This rotation matrix is a function of time and its evolution in time can be obtained as follows [12]:

$$
^I{\dot{\textbf{R}}}_B = {^I\textbf{R}_B{S({^B{\pmb{\omega}}_{B,I}})}}, \eqno{(1)}
$$
where $S({^B{\pmb{\omega}}_{B,I}})$ is the skew-symmetric matrix of angular velocity of the body with respect to the inertial frame as expressed in the body frame $^B{\pmb{\omega}}_{B,I}=[p,q,r]^T$. 

Likewise, the orientation of each motor frame $M_i$ can be obtained with respect to the body frame. First the position of the origin of frame $M_i$ with respect to body frame from the origin of the body frame can be written as:

$$
    ^B\textbf{O}_{M_i} = \textbf{R}_{z_B}((i-1)\frac{\pi}{2}) 
    \begin{bmatrix}
    L \\
    0 \\
    d \\
    \end{bmatrix} , (i=1,2,...,4), \eqno{(2)}
$$

Since we are using a quadcopter in ``+" configuration, we assume that the motors are evenly distributed by angle $\frac{\pi}{2}$ about axis $z_B$. Finally, the transformation from frame $M_i$ to body frame is obtained as follows:

$$
    ^B\textbf{R}_{M_i} = \textbf{R}_z((i-1)\frac{\pi}{2})\textbf{R}_y(\beta_i)\textbf{R}_x(\alpha_i), (i=1,2,...,4), \eqno{(3)}
$$

\subsection{Equations of Motion}

The quadcopter is consisted of several rigid bodies and it is considered to be symmetric about its axes of rotation. Because of the symmetry, the inertia tensor of the vehicle, $I^B$, will be diagonal and is expressed in the body frame. We also assume that the moment of inertia of the propellers, $I^p$, are very small compared to $I^B$. We can neglect drag force in angular motion of the body by assuming very small angular velocities. Considering these simplifying assumptions, the rotational motion is governed by the following equation: 

$$
    \pmb{\tau}=I{^B\dot{\pmb{\omega}}_{B,I}}+{^B\pmb{\omega}_{B,I}}\times{(I{^B\pmb{\omega}_{B,I}}+\sum_{i=1}^{4}{I^p{\pmb{\omega}^{p_i}}})}, \eqno{(4)}
$$
where $\pmb{\omega}^{p_i}$ is the angular velocity of the propeller with respect to the inertial frame as expressed in the body frame. $\pmb{\tau}$ is the torque generated by thrust forces and the reaction from motors expressed in body frame. Thrust force and reaction torque of each propeller $P_i$ in the frame $M_i$, can be approximated by the following formulas [13]: 

$$
    ^{M_i}\textbf{F}_{P_i}=[0,0,k_f{\dot{\gamma}_i^2}]^T, \eqno{(5)}
$$
$$
    ^{M_i}\pmb{\tau}_{P_i}=(-1)^{i+1}k_t{^{M_i}\textbf{F}_{P_i}}, \eqno{(6)}
$$
Using (5) and (6), we have: 

$$
    \pmb{\tau} = \sum_{i=1}^{4}({^B\textbf{O}_{M_i}}\times{{^B\textbf{R}_{M_i}}{^{M_i}\textbf{F}_{P_i}}}+{^B\textbf{R}_{M_i}}{^{M_i}\pmb{\tau}_{P_i}}), \eqno{(7)}
$$

The position of the vehicle in inertial frame is shown by Cartesian coordinates \textbf{s}$=[s_1,s_2,s_3]^T$. Finally, the equation governing translational motion can be written as follows: 

$$
    m\ddot{\textbf{s}} = {^I\textbf{R}_B}\sum_{i=1}^{4}({{^BR_{M_i}}{^{M_i}\textbf{F}_{P_i}}})+m\textbf{g}, \eqno{(8)}
$$
where $m$ is total mass of the vehicle and \textbf{g} is gravitational acceleration vector expressed in the inertial frame. 

\section{EFFECTS OF DIHEDRAL AND TWIST ANGLES}

In this part, an aerodynamic phenomenon, called dihedral effect, which is very common in fixed wing aircraft is introduced [14]. As shown in Fig. 4, when quadcopter is hovering, local air linear velocity with respect to the blade is equal to $\dot{\gamma_i}{C_{blade}(r)}$, where ${C_{blade}(r)}$ is the distance from the blade element (airfoil) to the shaft of the motor. At hover condition, it is assumed, this is the only relative velocity between the blade and the air. In this case, the angle of attack (AOA) of the blade, $\Theta_i$, will be defined as the angle between the chordline of the blade element and the velocity vector of the airflow over the blade that is shown in blue color in Fig. 4.

   \begin{figure}[b]
      \centering
      \includegraphics[scale=0.8]{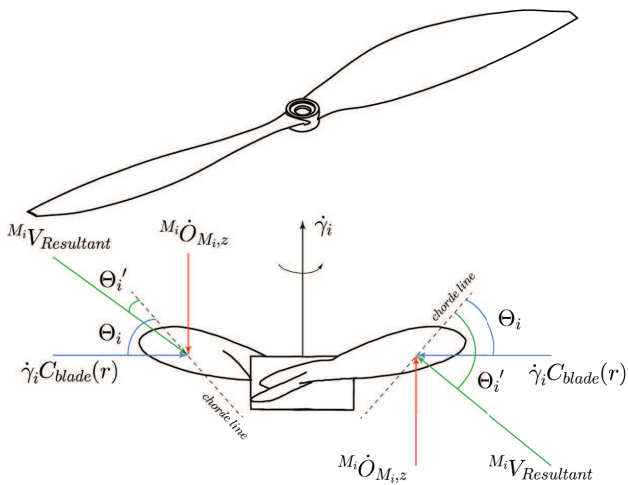}
      \caption{Dihedral Effect - On top is a propeller and on bottom is a front view of it. In the left, is the case when moving the motor up and in the right, is the case when moving the motor down.}
      \label{betai}
   \end{figure}

Moving the motor up or down, will generate an additional relative velocity between the blade and airflow, which in this case is parallel to the angular velocity of the propeller and is shown in red color in Fig. 4. It should be noted that in this figure, for visualization purpose and to save space, if the motor is moving down, the dihedral effect is shown in the right side of the figure and if the motor is moving up, dihedral effect is shown in the left side of the figure. The resultant of this additional velocity of airflow (due to the translational movement of the motor as shown in red in Fig. 4) with the linear velocity of airflow at each blade element due to rotation of the propeller (as shown in blue in Fig. 4) is the total airflow velocity relative to the blade, $^{M_i}V_{Resultant}$ (shown in green in Fig. 4). 

If the motor is moving down (see right side of Fig. 4), it is clear that the AOA increases and as a result, thrust force $^{M_i}F_{P_i}$ will increase [14]. On the other hand, if the motor is moving up (see left side of Fig. 4) the resultant velocity makes a smaller AOA than before when the motor was at rest, and as a result, thrust force decreases. This effect is called ``Dihedral Effect". In summary: 
\begin{itemize}
\item Any air flow with positive (negative) z-component velocity in frame $M_i$ increases (decreases) the AOA which increases (decreases) thrust force.
\end{itemize}

Now consider a quadcopter in a 2D motion. In Fig. 5, a configuration with no twist angle ($\alpha_i=0$) and constant dihedral angle $\beta_i=b$ ($b$ is negative) is shown. Consider the vehicle is pitching down and moving to the left which is equivalent of having an air flow with horizontal velocity to the right as shown in blue color in Fig. 5. According to ``Dihedral Effect", for the left motor, there will be an airflow with positive z-component in the frame $M_i$ as shown in green color and similarly, in the right motor, there will be an airflow with negative z-component velocity in the corresponding frame $M_i$. As a result, the AOA in the left motor increases thus its thrust force increases. But in the right motor, the AOA decreases and thrust force decreases as well. This interesting effect can make the vehicle stable in translational motion. As the vehicle moves to the left, due to the difference between thrust of the left and right motors, a moment $q'$, is generated that acts like damping in the system which resists with pitch motion and tries to reduce pitch angle to zero. 

   \begin{figure}[b]
      \centering
      \includegraphics[scale=0.4]{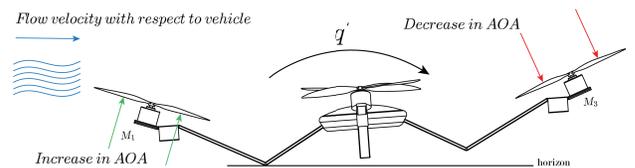}
      \caption{Dihedral effect in 2D motion of quadcopter. The quadcopter is pitching down and moving to the left. Dihedral effect generates the moment $q'$ and acts like a damping the in system. }
      \label{dihed_quad}
   \end{figure}
   
In order to derive equations for this force and moment, first we find the equation to calculate the thrust force as a function of AOA as follows [14]:
$$
    ^{M_i}\textbf{F}_{P_i} = {\frac{1}{2}}{\int_0^{C_{blade}}}{{\rho}v^2(r){C_l(r)}{c(r)}dr}, \eqno{(9)}
$$
where $\rho$ is the air density, $v(r)$ is the linear velocity of the airflow due to rotation of the blade at distance $r$ from the motor shaft, $C_l(r)$ is the lift coefficient of the blade element at distance $r$ from the motor shaft and $c(r)$ is the chord of blade element. 

At low speed flight, $C_l$ changes linearly with AOA [14], which can be written as: 
$$
    \frac{\Delta{C_l}}{\Delta{\Theta_i}} = \sigma, \eqno{10)}
$$

The linear velocity of $i^{th}$ motor with respect to the inertial frame as expressed in frame $M_i$ can be written as: 
$$
	^{M_i}{\dot{\textbf{O}}}_{M_i,I} = {^{M_i}[{\dot{O}}_{M_i,x},{\dot{O}}_{M_i,y},{\dot{O}}_{M_i,z}]^T}, \eqno{(11)}
$$

If $\dot{O}_{M_i,z}$ in (11) is positive, AOA and thrust force will decrease and if it is negative, we will have an increase in AOA and thrust force accordingly. Using Fig. 5 and trigonometric relations, we can find the change in AOA as follows: 
$$
	\Delta{\Theta_i} = \Theta_i - \Theta_i' = arctan(\frac{{\dot{O}}_{M_i,z}}{{\dot{\gamma_i}}{r}}), \eqno{(12)}
$$

Combining (9)-(12), and assuming hover conditions, ${\dot{O}}_{M_i,z}\ll{|\dot{\gamma_i}|r}$ and constant chord in the blade will result in the following: 
$$
	^{M_i} {\Delta}\textbf{F}_{P_i} = [0,0,-{\frac{1}{4}}c{\sigma}{\rho}{\dot{O}}_{M_i,z}{|{\dot{\gamma_i}|}{C_{blade}^2}}]^T, \eqno{(13)}
$$
where the negative sign in the equation together with the sign of $\dot{O}_{M_i,z}$ determine either the change in thrust is negative or positive. Near hover conditions, if we consider $\dot{\gamma_i}$ to be constant then (13) can be simplified further as follows:
$$
	^{M_i}{\Delta}\textbf{F}_{P_i} = [0,0,-\zeta\dot{O}_{M_i,z}]^T, \eqno{(14)}
$$
where $\zeta$ is a constant and a function of physical parameters of the blade and the airflow. Equation (14) is called ``pitch damper" and likewise we will have a ``roll damper". 

Effects of $\alpha_i$ is also in the category of ``Dihedral Effects". As shown in Fig. 6, to damp yaw motion, we need to choose $\alpha_{1,3}>0$ and $\alpha_{2,4}<0$. This is an interesting case where dihedral effect damps yaw motion. To better visualize this effect, assume that the quadcopter has a positive rotation about axis $z_B$. In this case, due to dihedral effect, AOA of propellers 2 and 4 decreases since there is an airflow with negative z-component of its linear velocity in the corresponding frames $M_{2,4}$. On the other hand, there will be an air flow with positive z-component of its linear velocity in the corresponding frames $M_{1,3}$. This phenomenon, generates a yaw moment (shown in green) that resists yaw motion $r$. Note that using $\alpha_{1,3}<0$ and $\alpha_{2,4}>0$ will have an adverse effect on yaw motion and could make it unstable. 
   \begin{figure}[b]
      \centering
      \includegraphics[scale=0.45]{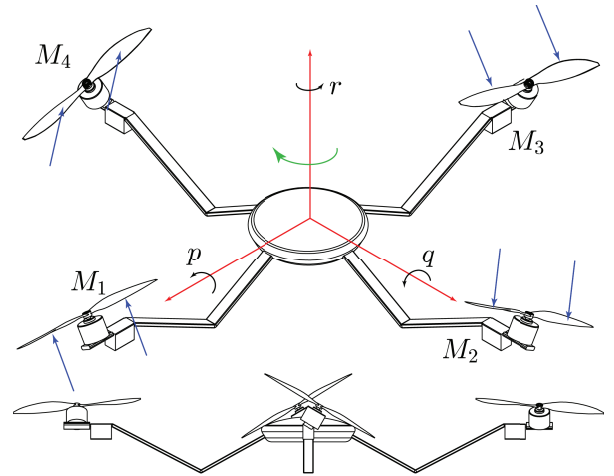}
      \caption{Quadcopter having only twist angles $\alpha_{1,3}>0$ and $\alpha_{2,4}<0$. The vehicle is going through pure yaw motion $r$ and dihedral effect generates a counteracting yaw motion that damps yaw motion.}
      \label{alpha1}
   \end{figure}

Finally, we will have changes in motor thrusts according to the following equations: 
$$
	^{M_i}{\Delta}\textbf{F}_{P_i,roll} = [0,0,-\zeta_{roll}\dot{O}_{M_i,z}]^T, \eqno{(15)}
$$
$$
	^{M_i}{\Delta}\textbf{F}_{P_i,pitch} = [0,0,-\zeta_{pitch}\dot{O}_{M_i,z}]^T, \eqno{(16)}
$$
$$
	^{M_i}{\Delta}\textbf{F}_{P_i,yaw} = [0,0,-\zeta_{yaw}\dot{O}_{M_i,z}]^T, \eqno{(17)}
$$
$$
	^{M_i}{\Delta}\textbf{F}_{P_i} = {^{M_i}{\Delta}\textbf{F}_{P_i,roll}}+{^{M_i}{\Delta}\textbf{F}_{P_i,pitch}}+{^{M_i}{\Delta}\textbf{F}_{P_i,yaw}}, \eqno{(18)}
$$

These changes in thrust force of each motor will affect translational motion as well as rotational motion by generating a moment about COM of the vehicle.

\section{STABILITY ANALYSIS}

In this section, we expand the equations for yaw motion and present the effects of twist angle on stability in yaw motion followed by a discussion on effects of dihedral angle in roll and pitch motion and also the effects of location of center of mass on overall stability of the vehicle. Also, for brevity, cross-coupling of the angular momentum of propellers are not presented in this section. At the end six different configurations using aforementioned parameters will be compared in terms of stability and maneuverability. 

\subsection{Effect of Twist Angle in Yaw Motion}

For simplicity assume that $\beta_i=0$, $\alpha_{1,3}>0$ and $\alpha_{2,4}<0$. Also all of these angles are kept constant during the analysis.  Here $d$ is positive meaning that the center of mass is located below the flat plane. Using (4) and (7), the equations governing the rotational motion can be written as follows: 
$$
	\pmb{\tau} = 
    \begin{bmatrix}
    {I_{xx}\dot{p}} \\
    {I_{yy}\dot{q}} \\
    {I_{zz}\dot{r}}
    \end{bmatrix} +
    \begin{bmatrix}
    (I_{zz}-I_{yy})qr \\
    (I_{xx}-I_{zz})pr \\
    0
    \end{bmatrix}
    , \eqno{(19)}
$$
where
$$
	\pmb{\tau} = 
    \begin{bmatrix}
    k_fds_a({\dot{\gamma_1}^2} - {\dot{\gamma_3}^2})+(k_fLc_a+k_tk_fs_a)({\dot{\gamma_2}^2} - {\dot{\gamma_4}^2}) \\
    k_fLc_a({\dot{\gamma_4}^2} - {\dot{\gamma_2}^2})+(k_tk_fs_a+k_fds_a)({\dot{\gamma_3}^2} - {\dot{\gamma_1}^2}) \\
    (k_tk_fc_a-k_fLs_a)({\dot{\gamma_1}^2} - {\dot{\gamma_2}^2} + {\dot{\gamma_3}^2} - {\dot{\gamma_4}^2})
    \end{bmatrix}
$$
where $\tau$ is the external torque generated by the motors to control the attitude of the vehicle and $\dot{\gamma}_i$ is the RPM of the motors. Also $s$ and $c$ represent sine and cosine functions.  It can be shown that $\alpha_i=0$ yields equations of motion for a regular quadcopter without tilting angles (details are saved for the sake of brevity). From (19), in a pure yaw motion, we have the following equation: 
$$
	\tau_{yaw}=I_{zz}\dot{r}, \eqno{(20)}
$$

Assuming the motors input for yaw motion to be equal to $u={\dot{\gamma_1}^2} - {\dot{\gamma_2}^2} + {\dot{\gamma_3}^2} - {\dot{\gamma_4}^2}$, we can rewrite (19) as follows:
$$
	\dot{r}=\frac{(k_tk_fc_a-k_fLs_a)}{I_{zz}}u, \eqno{(21)}
$$

Taking Laplace transform of (21), we can derive yaw motion transfer function as follows:

$$
	\frac{r(s)}{u(s)}=\frac{C_1}{s}, \eqno{(22)}
$$
where $C_1=\frac{(k_tk_fc_a-k_fLs_a)}{I_{zz}}$. 

Using (17), we can add the effects of twist angle into (21). Suppose the vehicle is going through pure yaw motion, $r$, as shown in Fig. 6. This yaw motion will generate local airflow over each blade with linear velocity equal to:  
$$
	^B\textbf{v}_{P_i}=[0,0,r]^T\times{^B\textbf{O}_{M_i}}, \eqno{(23)}
$$

Using (17) and (23), one can calculate the change in thrust force for each motor: 
$$
	^{M_i}{\Delta}\textbf{F}_{P_i,twist} = -\zeta_{yaw}{^B\textbf{R}_{M_i}^T}{^B\textbf{v}_{P_i}}, \eqno{(24)}
$$

For all motors, torque due to (24) can be calculated as: 
$$
	\pmb{\tau}_{twist} = \sum_{i=0}^4{{^B\textbf{O}_{M_i}}\times{^B\textbf{R}_{M_i}}{^{M_i}{\Delta}\textbf{F}_{P_i,twist}}}, \eqno{(25)}
$$

As shown in Fig. 6, any yaw motion $r$, will generate an airflow with negative z-component of  its linear velocity in the frames $M_{2,4}$.  Likewise, it will generate an airflow with positive z-component of its linear velocity in the frames $M_{1,3}$. As a result, based on (25), a torque will be generated that counteracts with the yaw motion $r$. Considering the simplifying assumptions made earlier in this section and using (23) we can calculate (24) as follows: 
$$
	^{M_i}{\Delta}\textbf{F}_{P_i,twist} = [0,0,(-1)^{i+1}\zeta_{yaw}Ls_ar]^T, \eqno{(26)}
$$
Using (25)- (26), we can write: 
$$
	\pmb{\tau}_{twist} = [0,0,-4\zeta_{yaw}L^2s_a^2r]^T, \eqno{(27)}
$$

Now, using (20) and the third component of (27) namely $\tau_{twist,yaw}$, we can add the effects of twist angle into equation (20) as follows: 
$$
	\tau_{yaw}+\tau_{twist,yaw} = I_{zz}\dot{r}, \eqno{(28)}
$$
$$
	C_1u = \dot{r}+\frac{\zeta'_{yaw}}{I_{zz}}r, \eqno{(29)}
$$
where $\zeta'_{yaw}=4\zeta_{yaw}L^2s_a^2>0$. Taking Laplace transform of (29) and simplifying it will result in: 
$$
	\frac{r(s)}{u(s)}=\frac{C_1}{s+{\frac{\zeta'_{yaw}}{I_{zz}}}}, \eqno{(30)}
$$

Comparing (30) with (22) shows that the vehicle has become more stable indeed. Transfer function (30) shows that it only has one negative pole indicating asymptotic stability in yaw motion. In addition to stability, this configuration helps to yaw faster because of the twist angle. Using twist angle a component of thrust force can be used to generate yaw motion which can be larger and easier to generate compared to regular quadcopters which yaw using only reaction torques of the motors. Note that $\alpha_{1,3}<0$ and $\alpha_{2,4}>0$ will have adverse effect on stability and will destabilize the system. 

Similarly, it can be shown that such phenomenon exists in roll and pitch motion for negative values of dihedral angle $\beta_i$ and similar transfer functions can be derived (details are not provided for brevity). The effect of location of center of mass is hidden in the value of $\zeta'$ in roll and pitch motion. It can be shown that for $d>0$ , as $d$ increases, the location of the pole of the transfer function will move to the left in the complex plane and increases stability and decreases maneuverability. Similarly, as $d$ decreases (even for negative values), the location of the pole of the transfer function will move to the right in the complex plane and stability will be decreased and maneuverability will be increased. 

\subsection{Comparison of Six Specific Configurations Based on Dihedral and Twist Angles}

In this section, based on dihedral and twist angles, six different configurations are proposed followed by a comparison in terms of stability and maneuverability. A regular quadcopter with all motors' angles set to zero is considered as a reference for comparison. The sign of dihedral and twist angles for each motor determines degree of stability or maneuverability in each configuration. The following list, ranks these configurations from the most stable to the least stable (for simplicity, we assume that $d$ is positive for all configurations): 
\begin{enumerate}
\item $\beta_i<0$, $\alpha_{1,3}>0$ and $\alpha_{2,4}<0$
\item $\beta_i<0$, $\alpha_i=0$
\item $\beta_i=0$, $\alpha_{1,3}>0$ and $\alpha_{2,4}<0$
\item $\beta_i=0$, $\alpha_i=0$
\item $\beta_i=0$, $\alpha_{1,3}<0$ and $\alpha_{2,4}>0$
\item $\beta_i>0$, $\alpha_{1,3}<0$ and $\alpha_{2,4}>0$
\end{enumerate}

In configuration (1), dihedral and twist angles are in favor of the stability and three dampers for roll, pitch and yaw motion are active in the quadcopter and are helping to stabilize its rotational motion. In configuration (2), twist angles are all set to zero, meaning that no damping (due to twist angles) exist in yaw motion and only roll and pitch dampers are active which results in having a vehicle less stable compared to configuration (1). In Configuration (3), only yaw damper is active and in configuration (4) all dihedral and twist angles are set to zero representing a regular quadcopter without tilting angles of the motors. In configuration (5), twist angles have an adverse effect compared to what we had in configuration (1), meaning that twist angles in this configuration will destabilize yaw motion of the quadcopter.

Note that having an adverse effect on stability means that the poles of the transfer function will move rightward in the complex plane and in some cases the poles will possibly fall in the right half of complex plane. Finally, in configuration (6), all dihedral and twist angles are having adverse effect with regard to stability in the system. However, in configurations (5) and (6), the vehicle has the highest maneuverability compared to other configurations. In summary, depending on applications and the environment in which the quadcopter is operating, choosing the best configuration and optimized values for dihedral and twist angles will be a trade off between stability and maneuverability. 

   \begin{figure}[t]
      \centering
      \includegraphics[scale=0.46]{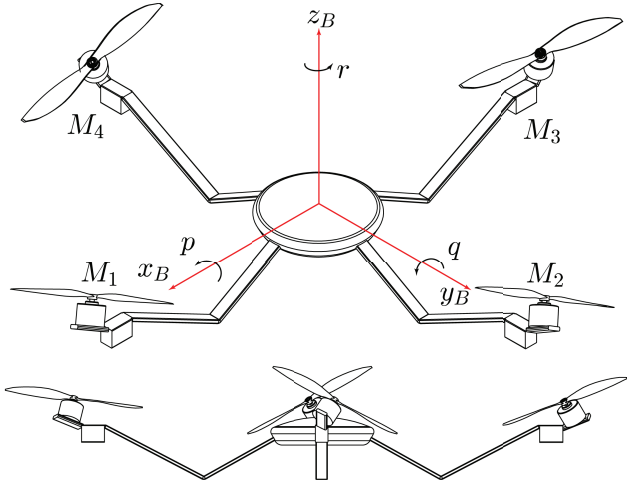}
      \caption{Quadcopter having both dihedral and twist angles. In this configuration $\beta_i<0$, $\alpha_{1,3}>0$ and $\alpha_{2,4}<0$ which renders the most stable configuration considering tilting angles of the motors. }
      \label{alpha}
   \end{figure}

\section{CONCLUSIONS}

Equations of motion of a quadcopter with tilted motors and having center of mass offset, in the z-direction of the body frame, were derived. The effects of tilting angles (dihedral and twist angles) on the thrust generated by propellers and consequently on stability of the system were introduced afterwards. Transfer functions considering pure yaw motion were derived followed by stability analysis and formulation of a yaw damper produced by adding twist angles to the motors for a specific configuration. Six different configurations based on these angles were introduced and were ranked based on stability and maneuverability. One of those configurations led to finding the most stable design (Fig. 7) with intrinsic damping in roll, pitch and yaw motion. The formulation for these dampers was presented followed by stability analysis in yaw motion. 

The dampers in the system would be favorable for applications where the vehicle hovers such as imaging, surveillance and monitoring. They will be unfavorable when the vehicle is in motion and maneuverability is needed. As seen in section IV, both stability and maneuverability can be achieved using different configurations. As a future work, a reconfigurable system can be designed in a way to transform from the most stable system to the most maneuverable system in the respective situation and vice versa. Such vehicle will be able to change dihedral and twist angles on the fly in order to transform to the required configuration. 

Another possible future work is to find the optimized values for dihedral and twist angles. Two different optimization problems can be defined: 1) optimizing the angles for the most stable configuration; and 2) optimizing the angles for the most maneuverable configuration. Finally, verifying the results of this paper using experiments will be done in a future work as well.

\addtolength{\textheight}{-12cm}   % This command serves to balance the column lengths

%%%%%%%%%%%%%%%%%%%%%%%%%%%%%%%%%%%%%%%%%%%%%%%%%%%%%%%%%%%%%%%%%%%%%%%%%%%%%%%%

\end{document}